\documentclass[b5paper,twoside,headrule,reqno,11pt,]{amsart}
\usepackage{latexsym, amsmath, amssymb, amsfonts, amscd}
\usepackage{graphics}
\usepackage[dvips]{epsfig}
\usepackage{latexsym}
\usepackage{amssymb}
\usepackage[flushmargin]{footmisc}

\newcommand{\gra}{\alpha}
\newcommand{\grb}{\beta}
\newcommand{\grd}{\delta}          
\newcommand{\gre}{\epsilon}
  
\newcommand{\grg}{\gamma}          \newcommand{\grG}{\Gamma}

\newcommand{\grk}{\kappa}
\newcommand{\grl}{\lambda}

\newcommand{\grr}{\rho}
\newcommand{\grs}{\sigma}          

\newcommand{\grw}{\omega}          \newcommand{\grW}{\Omega}
             
\newcommand{\gry}{\eta}
\newcommand{\grz}{\zeta}
\newcommand{\grth}{\theta}         

\newcommand{\grph}{\phi}

\newcommand{\del}{\partial}

                           \newcommand{\bC}{\mathbf{C}}

                           \newcommand{\bZ}{\mathbf{Z}}

\newcommand{\beq}{\begin{equation}}
\newcommand{\eeq}{\end{equation}}

\newtheorem{theorem}{Theorem}[section]
\newtheorem{lemma}[theorem]{Lemma}

\newcommand{\barr}{\overline}

\begin{document}


\def\evenhead{{\protect\centerline{\textsl{\large{Xianghong Gong and S. M. Webster}}}\hfill}}

\def\oddhead{{\protect\centerline{\textsl{\large{Regularity for The CR Vector Bundle Problem I}}}\hfill}}

\pagestyle{myheadings} \markboth{\evenhead}{\oddhead}

\thispagestyle{empty} \noindent{{\small\rm Pure and Applied
Mathematics Quarterly\\ Volume 6, Number 3\\ (\textit{Special Issue:
In honor  of \\ Joseph J. Kohn})\\
983---998, 2010} \vspace*{1.5cm} \normalsize

\begin{center}
\Large {\bf Regularity for The CR Vector Bundle Problem I}
\end{center}
\footnotetext{Received November 16, 2007.}
\renewcommand{\thefootnote}{\fnsymbol{footnote}}
\begin{center}
\large Xianghong Gong\footnote{Partially supported by NSF grant
 DMS-0705426.}  and S. M. Webster
\end{center}

\centerline{\it Dedicated to Professor J.J.Kohn on the occasion of
his 75th birthday}
\begin{center}
\begin{minipage}{5in}
\noindent {\bf Abstract:} We give a new solution to the local
integrability problem for CR vector bundles over strictly
pseudoconvex real hypersurfaces of dimension seven or greater.  It
is based on a KAM rapid convergence argument and avoids the previous
more difficult Nash-Moser methods.  The solution is sharp as to
H\"older continuity.
\end{minipage}
\end{center}

\tableofcontents

\vspace{3ex}


\setcounter{section}{0} \setcounter{equation}{0}

\noindent\textbf{Introduction}.  The idea of a CR vector bundle
$E\rightarrow M$, over a CR manifold $M$, naturally generalizes that
of a holomorphic vector bundle over a complex manifold.  Here we
shall be concerned with the local integrability problem over a
strongly pseudoconvex real hypersurface $M^{2n-1}\subset \bC^{n}$.
Therefore, we assume that $E$ is trivial of some rank $r$, i.\ e.\
$E\cong M\times\bC^{r}$, via a local frame field $e=(e_{1},\ldots
,e_{r})$, and we phrase the problem in terms of a connection $D$ on
$E$,
\begin{equation}
    De = \grw\otimes e, \; \; \; \; \grW = d\grw - \grw\wedge\grw.
\end{equation}

The integrability condition is that the curvature 2-forms $\grW =
(\grW_{i}^{\; j})$ belong to the ideal $\mathcal{J}(M)$ generated by
the restrictions of (1,0)-forms to $M$. The problem is to find a new
frame $e$, so that the connection 1-forms $\grw = (\grw_{i}^{\;j})$
also belong to $\mathcal{J}(M)$.  We focus mainly on the regularity
of the new frame.

Over complex manifolds the problem goes back to Koszul and Malgrange
[5]. Over CR manifolds it was considered in [13], in conjunction
with the Nash-Moser theory developed in [12] for the much more
difficult CR embedding problem. Ma and Michel [6], [7] have improved
the regularity for these results. Here we shall consider the vector
bundle problem in its own right.  We shall eliminate the difficult
Nash-Moser techniques and derive results sharp as to regularity via
a very natural KAM rapid convergence argument.

More precisely, we have $d\mathcal{J}(M)\subseteq\mathcal{J}(M)$,
and denote by $\barr{\del}_{b}$ the reduction of the exterior
derivative $d$ mod $\mathcal{J}(M)$. In section 1 we choose suitable
representatives $\barr{\del}_{M}$ for $\barr{\del}_{b}$, and
$\grw''$ for $\grw$, mod $\mathcal{J}(M)$. Since the other
components of $\grw$ will be irrelevant in this work, we simplify
the notation by setting $\grw=\grw''$. Then, for a frame change
$\tilde{e} = Ae$, $\det A\neq0$, the derivation property of $D$, and
reduction mod $\mathcal{J}(M)$ give
\begin{equation}
   \tilde{\grw}A = \barr{\del}_{M}A + A\grw.
\end{equation}
We want to make $\tilde{\grw} = 0$.  Thus, we seek a solution $A$ in
some neighborhood of any given point to
\begin{equation}
   - A^{-1}\barr{\del}_{M}A = \grw , \; \;\; \mbox{if}
   \; \; \; \barr{\del}_{M}\grw = \grw\wedge\grw,
\end{equation}
the latter expressing the integrability condition.  If the rank
$r=1$, this reduces to the $\barr{\del}_{b}$-problem, which we know
to be locally solvable, if $\dim M = 2n-1 \geq 5$, [4], [1].  If $r
> 1$, (0.3) casts the problem in a non-linear light, which reflects
our methods.

To measure regularity, we consider both the standard H\"older
spaces $C^{k,\gra}(M)$, and the Folland-Stein H\"older spaces
$C^{k,\gra}_{FS}(M)$ [2], for integer $k$ and $0 \leq \gra < 1$.
We assume that $M$ is of class $C^{l}$, $5\leq l \leq \infty$, and
that $\dim M = 2n-1 \geq 7$, $n\geq 4$  (thus, omitting $\dim M =
5$), and prove the following.
\begin{theorem}
 Suppose that $\grw$ is of class $C^{k}(M)$, $k_{1}\leq k \leq l-4$,
with $k_{1}=1$.  Then there exists a local solution $A$ of class
$C^{k,\gra}(M)$, for $0\leq\gra\leq 1/2$. In particular, if $M$ and
$\grw$ are of class $C^{\infty}$, there is a solution $A$ of class
$C^{\infty}$.
\end{theorem}
The proof given here is complete based on essentially known
estimates [4], [11].  Some improvement is still possible.  For
example, it will follow from [3] that, with appropriate weak
definitions, we can get a solution with $l = 3$ and $k = 0$.

Our second result is restricted to a real hyperquadric.
\begin{theorem}
  Let $\grw$ be of class $C^{k,\gra}_{FS}(M)$, $k_{2}\leq k < \infty$,
with $k_{2}=1$, and $0 < \gra < 1$, where $M$ is the Heisenberg
group. Then there exists a local solution $A$ of class
$C^{k+1,\gra}_{FS}(M)$.
\end{theorem}
This can be improved to $k_{2}=0$, and with some effort carried
over to vector bundles over general strictly pseudoconvex real
hypersurfaces. We shall pursue this in a future work.

As in [13], the main technical tool is the $\barr{\del}_{b}$ local
homotopy formula of Henkin [4],
\begin{equation}
  \grph^{0,q} = \barr{\del}_{M}P\grph + Q\barr{\del}_{M}\grph, \; \;
  0 < q < n-2,
\end{equation}
on suitable subdomains $M_{\grr}\subset M$.  We need it for $q=1$ as
in [12], which is why the $5$-dimensional case is omitted.  We also
refer to the result of Nagel and Rosay [9].  While a linear argument
is conceivable, our estimates for the operators $P$, $Q$ are not
adequate (see also Romero [10] in this respect).  Therefore, we
proceed (somewhat imprecisely) as follows.

We take $A=I+B$, $||B|| < 1$, apply (0.4) with $\grph = \grw$, and
use the integrability condition (0.3).  This gives
\begin{equation}
  \tilde{\grw}A =\barr{\del}_{M}(B + P\grw) + Q(\grw\wedge\grw) + B\grw.
\end{equation}
The choice $B=-P\grw$ gives formally $||\tilde{\grw}||\leq
c||\grw||^{2}$. This will be repeated on a sequence of shrinking
domains $M_{\grr}\subset M$. The rapid (quadratic) decrease of
$||\grw||$ will allow us to overcome deficiencies (loss of
derivatives and blow-up of coefficients in our estimates) and prove
convergence to a solution of the problem.

In section 2 we summarize some previous estimates, which are not
sharp but allow us to get directly to the core argument, which is
given in section 3.  To get the procedure started, we need
$||\grw||$ to be sufficiently small, initially.  This is achieved
via a Taylor polynomial argument in section 1. At first we see with
this an apparent loss of over half the derivatives in the argument.
However this is greatly reduced in section 4, where we show
following Moser [8] that, with rapid convergence in some low norm,
the higher order derivatives are automatically pulled in. This
gives a weaker form of theorem (0.1), where $A\in C^{k}$.  The full
theorem follows with the H\"older estimates of section 5.

Another way to avoid this derivative loss on initial shrinking is to
use scale invariant Folland-Stein norms, which we do in section 6 on
the Heisenberg group.  This gives theorem (0.2). The use of scaling
also gives another proof of theorem (0.1) on the Heisenberg group.

We point out that the main deficiencies here are the rather
imprecise estimates we have for the solution operators $P$, $Q$ on
shrinking domains. This is somewhat addressed in [3].

We add some brief historical remarks.  The central result on the
local integrability of structures is, of course, the
Newlander-Nirenberg theorem (1957).  Its importance is reflected in
the number of proofs it has since received.  The local embedding
problem for strictly pseudoconvex CR structures $M$ was first posed
by Kohn (1964).  Counter examples for dim$M = 3$ were given by
Nirenberg (1973).  The first positive results for dim$M\geq 7$ were
given by Kuranishi (1982), and Akahori (1987).

We re-emphasize that the the five-dimensional case is, at this time,
still unresolved for both the CR embedding problem and the CR vector
bundle problem.  Perhaps the case of a CR vector bundle over the
5-dimensional Heisenberg ball would be the easiest problem to
resolve.


\section{Initial normalization}
\setcounter{equation}{0}

Here we achieve the required initial smallness for the connection
form $\grw$, by  Taylor polynomial arguments.

We take our real hypersurface in the form
\begin{equation}
   M^{2n-1}\subset\bC^{n} : r(z) = -y^{n}+|z'|^{2} + h(z',x^{n})=0,
\end{equation}
where $z=(z',z^{n})$, $z^{n}=x^{n}+iy^{n}$, and
$h=O(|(z',x^{n})|^{3})$ is of class $C^{l}$, $l\geq 3$.  For a basis
of complex vector fields we take
\begin{equation}
  X_{\gra}= \del_{\gra}-(r_{\gra}/r_{n})\del_{n}, \; \; \;
 X_{\barr{\gra}}=\barr{X_{\gra}}, \; \; \; T = \barr{T}=\del_{x^{n}},
\end{equation}
where $r_{j}=\del_{j}r = \del r/\del z^{j}$, and small Greek indices
have the range $1\leq\gra,\grb\leq n-1$. The dual basis is
\begin{equation}
    dz^{\gra}, \; \; d\barr{z}^{\gra}, \; \; \grth = \barr{\grth}= -i\del r,
\end{equation}
and $\mathcal{J}(M)$ is the differential ideal generated by
$dz^{\gra},\grth$. With a standard multi-index notation
$A=(\gra_{1},\ldots,\gra_{q})$, we have representatives for a
$(0,q)$-form $\grph$ and $\barr{\del}_{b}\grph$,
\begin{equation}
   \grph^{q}=\grph^{(0,q)} = \sum_{|A|=q}\grph_{\barr{A}}d\barr{z}^{A},\; \;
  \;   \barr{\del}_{M}\grph = \sum
       \barr{\del}_{M}\grph_{\barr{A}}\wedge d\barr{z}^{A},
\end{equation}
and for function $f$,
\begin{equation}
  \barr{\del}_{M}f = \sum_{\gra=1}^{n-1} X_{\barr{\gra}}f d\barr{z}^{\gra}.
\end{equation}

For (the (0,1)-part of) the connection one-forms we put
\begin{equation}
   \grw = \sum_{\gra=1}^{n-1} \grG_{\barr{\gra}}d\barr{z}^{\gra}.
\end{equation}
Thus, the $\grG_{\barr{\gra}}$ are $r\times r$ matrices of
functions. The integrability condition (0.3) is equivalent to
\begin{equation}
  X_{\barr{\gra}}\grG_{\barr{\grb}} - X_{\barr{\grb}}\grG_{\barr{\gra}} =
  \grG_{\barr{\gra}}\grG_{\barr{\grb}} - \grG_{\barr{\grb}}\grG_{\barr{\gra}},
\end{equation}
while the change of frame formula (0.2) is
\begin{equation}
   \tilde{\grG}_{\barr{\gra}}A = X_{\barr{\gra}}A + A\grG_{\barr{\gra}}.
\end{equation}

If $\grw$ is of class $C^{k}$, $k\leq l-1$, we can arrange that it
vanish to order $k$ at $0$, by a standard Taylor polynomial
argument, as follows.  Assume that $\grw$ starts with terms of order
$s$,
\begin{equation}
  \grG_{\barr{\gra}} = \grG_{\barr{\gra}}^{(s)} + \cdots , \; \;
  \grG_{\barr{\gra}}^{(s)} =
  \sum_{0\leq|B|\leq s}\grG_{\barr{\gra},\barr{B}}(z',x^{n})\barr{z'}^{B},
\end{equation}
where $\grG_{\barr{\gra},\barr{B}}$ is homogeneous of degree
$s-|B|$, and symmetric in the indices $\barr{B}$. We make the change
(1.8) with
\begin{equation}
  A = I + A^{(s+1)}, \; \; A^{(s+1)} = \sum_{0\leq|B|\leq s+1}
  A_{\barr{B}}(z',x^{n})\barr{z'}^{B},
\end{equation}
where the coefficients are symmetric in $\barr{B}$. We need to make
\begin{equation}
   \del_{\barr{\gra}}A + \grG_{\barr{\gra}}^{(s)} = 0.
\end{equation}
This defines the $A_{\barr{\gra B}}$ in terms of the
$\grG_{\barr{\gra},\barr{B}}$, consistently, since
$\grG_{\barr{\gra},\barr{B}}$ is symmetric in all its indices by the
integrability condition (1.7).  Thus we have the following.
\begin{lemma}
  If the connection form $\grw$ is of class $C^{k}$, $k\leq l-1$, then
we can arrange $\grw = O(k)$ at $0\in M$, by a preliminary change of
frame.
\end{lemma}

We have an analogous result relative to Folland-Stein derivatives on
the real hyperquadric (Heisenberg group), where $h(z',x^{n})=0$.  We
say that a function $f$ is of class $C^{k}_{FS}$ if it has
continuous derivatives
\begin{equation}
   T^{m}X^{S}\barr{X}^{R}f, \; \; 2m + |S| + |R| \leq k,
\end{equation}
where $X^{S}= X_{1}^{s_{1}}\cdots X_{n-1}^{s_{n-1}}$, etc.
\begin{lemma}
  Suppose that $\grw$ is of class $C^{k}_{FS}$ on the Heisenberg group $M$.
Then by a preliminary change of frame, we can achieve
\begin{equation}
 T^{m}X^{S}\barr{X}^{R}\grG_{\barr{\gra}}(0) = 0, \; \; 2m + |S| + |R| \leq k.
\end{equation}
\end{lemma}
The proof is similar.  We assume that (1.13) holds for $s-1$, $s\leq
k$, in place of $k$,
\begin{equation}
   \grG_{\barr{\gra}} = \grG_{\barr{\gra}}^{(s)} + \cdots ,
\end{equation}
where now $(s)$ refers to weight,
\begin{equation}
   wt(\barr{z'}^{R}z'^{S}(z^{n})^{m}) = |R| + |S| + 2m.
\end{equation}
We then take $A= I + A^{(s+1)}$ with
\begin{equation}
   T^{m}X^{S}\barr{X}^{R}X_{\barr{\gra}}A^{(s+1)}(0) =
  - T^{m}X^{S}\barr{X}^{R}X_{\barr{\gra}}\grG_{\barr{\gra}}^{(s)}(0).
\end{equation}
The integrability condition (1.7) guarantees symmetry in the barred
coefficients, so that this is consistent.  Continuing this, we
achieve the lemma.


\section{Estimates for the homotopy formula}
\setcounter{equation}{0}

Here we give estimates for the operators $P$, $Q$ which can be
readily derived from known results.  We refer to Henkin [4], and
many references therein, as well as to [11].  For convenience we
follow the notation of [11].  This will allow us to get quickly to
the main convergence argument, which is given in the next section.

We work with the real hypersurfaces
\begin{equation}
  M_{\grr} = M\cap\{(x^{n})^{2} + y^{n} \leq \grr^{2}\}.
\end{equation}
For $\grr > 0$ sufficiently small, this is a graph over
$D_{\grr}\subset\mathbf{R}^{2n-1}$, which is approximately a ball of
radius $\grr$.  We consider two such domains
$D_{0}=D_{\rho_0}\subset\subset D=D_\rho$ with $\grr_{0} =
\grr(1-\grs)$, $0 < \grs < 1$.  Let $\grd(D)\approx\grr$ denote the
diameter, and $\grd(D_{0},\del D)\approx\grr\grs$ distance to the
boundary.  In the notation of [11], we have $P = P_{0} + P_{1}$ and
$Q = Q_{0} + Q_{1}$, where $P_{0}$, $Q_{0}$ are integral operators
over $M_{\grr}$, and $P_{1}$, $Q_{1}$ are integral operators over
the boundary $\del M_{\grr}$.

For $\grph\in C_{0}^{k}(D)$, we follow the procedure of sections 3
and 4 of [11].  In particular, we have (3.7) of [11], but without
the boundary integral.  Any $k$-th order derivative $\del^{k}$ on
$D$ can be expressed in terms of the vector fields (1.39) of [11].
We may summarize the procedure by writing
\begin{equation}
  \del^{k} P_{0}\grph = \sum_{|K|\leq k} P_{0}^{(K)}(\del^{K}\grph).
\end{equation}
Here the right hand side is a sum of certain partial derivatives of
order up to $k$, and the $P_{0}^{(K)}$ are operators of type not
worse than $P_{0}$, with kernels involving at most $k+3$ derivatives
of $r$.  The absolute integrability of these kernels follows from
[4].   This gives
\begin{equation}
  ||P_{0}\grph||_{C^{k}(D)} \leq c_{k+3}\grd(D)\|\grph\|_{C^{k}(D)},
\end{equation}
where $c_{k+3}$ is a constant depending on $k+3$ derivatives of $r$.
We combine this estimate with a cut-off function $\grl_{0} +
\grl_{1} = 1$ on $\barr{D}$, $\grl_{0}=1$ on $D_{0}$, $\grl_{0}=0$
near $\del D$. For $\grph\in C^{k}(\barr{D})$, the product rule and
$|\del\grl_{0}|\leq c/\grd(D_{0},\del D)$ give
\begin{equation}
  \|\grl_{0}\grph\|_{C^{k}(D)} \leq
  c \sum_{j=0}^{k}\grd(D_{0},\del D)^{j-k}\|\grph\|_{C^{j}(D)}.
\end{equation}

To estimate $P_{0}(\grl_{1}\grph)$ on $D_{0}$ we (crudely) let the
derivatives fall on the kernel, increasing the blow-up at the
boundary;
\begin{equation}
  \|P_{0}(\grl_{1}\grph)\|_{C^{k}(D_{0})} \leq
   c_{k+3}\frac{vol(D-D_{0})}{\grd(D_{0},\del D)^{2n+2k-1}}
                  \|\grph\|_{C^{0}(D)}.
\end{equation}
We estimate the boundary integral $P_{1}$ similarly as in (3.8),
(4.6), and (4.7) of [11], using such bounds as
$|r_{\grz}\cdot(\grz-z)|\geq c|\grz-z|^{2}$, $|r_{\grz}-r_{z}|\leq
c|\grz-z|$, $|\zeta^n-z^n|\geq c\rho^2\sigma$, and $d\zeta^n\wedge
r_z\cdot d\zeta\wedge r_{\zeta}\cdot d\zeta=d\zeta^n\wedge
r_{z'}\cdot d\zeta'\wedge r_{\zeta'}\cdot d\zeta'=O(\rho^2)$,
 for terms appearing in the kernel (1.31) of [11].
This gives
\begin{equation}
  \|P_{1}(\grph)\|_{C^{k}(D_{0})} \leq
  c_{k+3}\frac{vol(\del D)}{\sigma\grd(D_{0},\del D)^{2n+2k-2}}\|\grph\|_{C^{0}(D)}.
\end{equation}
We note that $vol(\del D)\leq c\grr^{2n-2}$ and $vol(D-D_{0})\leq
c(\grr\grs)\grr^{2n-2}$. Combining and simplifying gives
($\|\cdot\|_{\grr,k} = \|\cdot\|_{C^{k}(D_{\grr})}$)
\begin{equation}
  \|P(\grph)\|_{\grr(1-\grs),k} \leq c_{k+3}'
  \{\grr\sum_{s=0}^{k}(\grr\grs)^{s-k}\|\grph\|_{\grr,s} +
  \gry'\|\grph\|_{\grr,0}\}
  \leq \gry \|\grph\|_{\grr,k},
\end{equation}
where
\begin{equation}
  \gry = \gry(\grr,\grs,k) = c_{k+3}\gry', \; \;
     \gry' = \grs^{-2n-2k+1}\grr^{-2k}.
\end{equation}
Notice that for $k=0$, we don't need the cutoff function, nor (2.4),
(2.5) and the factors $\grr$ cancel in (2.6) giving (2.8) with $k =
0$.

A similar estimate holds for $Q$. These estimates on shrinking
domains will be somewhat improved in [3].


\section{A KAM rapid convergence argument}
\setcounter{equation}{0}

In this section we construct  a sequence of approximate solutions to
the problem of theorem (0.1) of the introduction and prove
convergence.  With the estimates of the previous section, we get a
solution with apparent derivative loss in $C^{k}$ spaces.

We define a sequence of \lq\lq radii\rq\rq, $\grr > 0$, by
\begin{equation}
  \grr_{j+1} = \grr_{j}(1 - \grs_{j}), \; \; \grs_{j}=2^{-j-1}, \; \; j\geq 0,
\end{equation}
which decrease to a positive limit $\grr_{\infty}$ depending on
$\grr_{0}>0$, which will be chosen later.  Then $M_{j}\equiv
M_{\grr_{j}}$, $0\leq j \leq\infty$, form a decreasing family of
neighborhoods of $0\in M$.

We consider frames $e_{j}$ for $E$ over $M_{j}$, and the associated
connection matrix $\grw_{j}$ of $(0,1)$-forms, $0\leq j < \infty$.
They will be chosen so that $e_{0}=e$ is the original frame,
normalized so that $\grw_{0}=\grw$ vanishes to suitably high order
at $0$, and for $j\geq 1$,
\begin{equation}
  e_{j} = A_{j}e_{j-1}, \; \; A_{j} = I + B_{j}, \; \;
  G_{l}= A_{l}A_{l-1}\cdots A_{1}.
\end{equation}
Then on $M_{l}$ we shall have
\begin{equation}
  \grw_{l}A_{l} = \barr{\del}_{M}A_{l} + A_{l}\grw_{l-1},
\end{equation}
and inductively,
\begin{equation}
  \grw_{l}G_{l} = \barr{\del}_{M}G_{l} + G_{l}\grw_{0}.
\end{equation}
We must show that we have convergence in $C^{k}(M_{\infty})$, $k\geq
1$,
\begin{equation}
  \grw_{l}\rightarrow 0, \; \; \; G_{l}\rightarrow G_{\infty},
\end{equation}
with the matrix of functions $G_{\infty}$ invertible.  Then
$A=G_{\infty}$ will be a solution to our problem.

On each $M_{j}$, $0\leq j < \infty$, we have the homotopy formula
(0.4) with operators $P_{j}$, $Q_{j}$.  In (0.5) we take $A_{j+1} =
I + B_{j+1}$, with
\begin{equation}
  B_{j+1} = -P_{j}\grw_{j},
\end{equation}
giving
\begin{equation}
  \grw_{j+1}A_{j+1} = Q_{j}(\grw_{j}\wedge\grw_{j}) + B_{j+1}\grw_{j}.
\end{equation}
From (3.6), (2.7) we get
\begin{equation}
    \|B_{j+1}\|_{\grr_{j+1},k} \leq \gry_{j}\|\grw_{j}\|_{\grr_{j},k},\; \; \;
 \;\; \gry_{j}\equiv\gry_{j}^{(k)} = c_{k+3}\grs_{j}^{-2n-2k+1}\grr_{j}^{-2k}.
\end{equation}

For $r\times r$ matrices of functions and a given $k$, there is a
constant $\tilde{c}_{k}\geq 1$ for which
\begin{equation}
  \|AB\|_{\grr,k}\leq \tilde{c}_{k}\|A\|_{\grr,k}\|B\|_{\grr,k},
\end{equation}
\begin{equation}
  \|(I+B)^{-1}\|_{\grr,k}\leq (1 - \tilde{c}_{k}\|B\|_{\grr,k})^{-1},
\end{equation}
provided $\tilde{c}_{k}\|B\|_{\grr,k} < 1$. (Note that (3.9) holds
when $A,B$ are $r\times r$ matrices of forms.) Thus, we shall want
to arrange
\begin{equation}
     \tilde{c}_{k}\|B_{j+1}\|_{\grr_{j+1},k} < 1/2, \; \;
     \|A_{j+1}^{-1}\|_{\grr_{j+1},k} \leq 2 ,
\end{equation}
for $0\leq j <\infty$, at least for $k=0$, to carry out the
procedure.

Given (3.11), we get from (3.7), (2.7)
\begin{eqnarray}
  \|\grw_{j+1}\|_{\grr_{j+1},k} & \leq & 2\{
  \|Q_{j}(\grw_{j}\wedge\grw_{j})\|_{\grr_{j+1},k} +
  \tilde{c}_{k}\|B_{j+1}\|_{\grr_{j+1},k}\|\grw_{j}\|_{\grr_{j+1},k} \} \\
    & \leq & \gry_{j}\|\grw_{j}\|_{\grr_{j},k}^{2}, \nonumber
\end{eqnarray}
where we have absorbed a factor of $4\tilde{c}_{k}$ into the
constant $c_{k+3}$.

For $j=0$ we assume that $\grw_{0}$ is of class $C^{m}$, $m > k$,
and normalize as in section 1 so that $\grw_{0}$ and all its
derivatives of order $m$ and less vanish at $0$.  Then
\begin{equation}
  \|\grw_{0}\|_{\grr_{0},k} \leq c\grr_{0}^{m-k}, \; \; \;
  \gry_{0}\|\grw_{0}\|_{\grr_{0},k} \leq c\grr_{0}^{m-3k},
\end{equation}
provided $m > 3k$.  Then by shrinking $\grr_{0}$ we get (3.11) for
$j = 0$.

We set
\begin{equation}
  \grd_{j} \equiv \grd_{j}^{(k)} = \|\grw_{j}\|_{\grr_{j},k},  \; \; \;
  \grd_{j+1} \leq \gry_{j}\grd_{j}^{2}.
\end{equation}
We readily see that
\begin{equation}
  \gry_{j+1} = \gra_{j}\gry_{j}, \; \; \;
  \gra_{j} \equiv \gra_{j}^{(k)} = 2^{2n+2k-1}(1 - \grs_{j})^{-2k}.
\end{equation}
Since the $\grs_{j}$ decrease to $0$, the $\gra_{j}$ decrease to
$2^{2n+2k-1}$.  The $\gry_{j}$ increase to infinity. Finally we
define $\grz_{j}$,
\begin{equation}
   \grz_{j} \equiv \grz_{j}^{(k)} = \gra_{j}\gry_{j}\grd_{j}, \; \; \;
   \grz_{j+1}\leq \grz_{j}^{2}.
\end{equation}
By shrinking $\grr_{0}$, we arrange $\grz_{0}<1/2\tilde{c}_{k}$,
then the $\grz_{j}$ decrease rapidly to $0$, as also do   the
$\grd_{j}$ and $\|B_{j}\|_{\grr_{j},k}$, and we have (3.11) for all
$j\geq 0$.

From (3.11) we clearly have
\begin{equation}
  \|A_{j}\|_{\grr_{j},k} \leq 2, \; \; \;
  \|G_{l}\|_{\grr_{l},k}\leq \tilde{c}_{k}^{l-1}2^{l}.
\end{equation}
Since $G_{l}-G_{l-1} = B_{l}G_{l-1}$, we have
\begin{equation}
  \|G_{l}-G_{l-1}\|_{\grr_{l},k} \leq
  (2\tilde{c}_{k})^{l-1}\|B_{l}\|_{\grr_{l},k}
    \leq (2\tilde{c}_{k})^{l-1}\grz_{l-1}.
\end{equation}
It follows by the ratio test that
\begin{equation}
G_{l}-G_{1} = \sum_{j=2}^{l}(G_{j}-G_{j-1})
\end{equation}
converges in $C^{k}(M_{\infty})$ to a limit $G_{\infty}\in
C^{k}(M_{\infty})$. Furthermore, for any $\gre > 0$ we may arrange
$\grz_{0}<\gre/2\tilde{c}_{k}$, by shrinking $\grr_{0}$ as above.
Then $\grz_{l} < (\gre/2\tilde{c}_{k})^{l}$, and
\begin{equation}
    \sum_{l=1}^{\infty} (2\tilde{c}_{k})^{l+1}\grz_{l} \leq
       2\tilde{c}_{k}\gre/(1-\gre).
\end{equation}
Thus, by shrinking $\grr_{0}$ a second and final time, we ensure
that $G_{\infty}-G_{1} = G_{\infty}- I - B_{1}$ and $B_{1}$ are so
small that $G_{\infty}$ is invertible.

Thus, we have achieved a solution $A=G_{\infty}$ of class $C^{k}$,
in theorem(0.1) provided $\grw$ is of class $C^{3k+1}$.  This extra
smoothness requirement, which was needed only in (3.16) for the
initial smallness, will be removed in the following sections.


\section{Higher order derivatives}
\setcounter{equation}{0}  Now we assume that our initial
$\grw=\grw_{0}$ is of class $C^{k}$, $1\leq k\leq l-3$.  We
normalize so that it and all its first order derivatives vanish at
$0$.  By shrinking the initial radius $\grr_{0}$, we can make
$\grz_{0}^{(s)}$ arbitrarily small for $s=0$.  This allows us to
carry out the preceding argument.  The constructed sequences
$B_{j}$, $\grw_{j}$ are of class $C^{k}$, but are converging rapidly
to $0$ only in $C^{0}$-norm, and we cannot yet pass to the limit in
(3.4). By modifying an idea in Moser [8], we show that,
\emph{without any further change of the sequence}, we have
convergence $\grw_{l}\rightarrow 0$, $G_{l}\rightarrow G_{\infty}$
in $C^{k}(M_{\infty})$.

For this we take any first order partial derivative $\del_{x}$ in
(3.7),
\begin{equation}
   \del_{x}\grw_{j+1}A_{j+1} =
   \del_{x}Q_{j}(\grw_{j}\wedge\grw_{j}) +
   \del_{x}B_{j+1}\grw_{j} + B_{j+1}\del_{x}\grw_{j} -
   \grw_{j+1}\del_{x}B_{j+1}.
\end{equation}
We multiply this by $A_{j+1}^{-1}$ and take
($\|\cdot\|_{\grr_{j+1},0}$)-norms over $M_{j+1}$, using
\begin{eqnarray}
   \|\grw_{j+1}\|_{\grr_{j+1},0}  & \leq  & \gry_{j}^{(0)}
   \|\grw_{j}\|_{\grr_{j},0}^{2} , \\
   \|B_{j+1}\|_{\grr_{j+1},0}  & \leq  & \gry_{j}^{(0)}
   \|\grw_{j}\|_{\grr_{j},0}, \; \;
   \|A_{j+1}^{-1}\|_{\grr_{j+1},0}  \leq 2, \nonumber \\
   \|\del_{x}B_{j+1}\|_{\grr_{j+1},0}  & \leq  & \gry_{j}^{(1)}
   \|\grw_{j}\|_{\grr_{j},1} ,  \nonumber \\
   \|\del_{x}Q_{j}(\grw_{j}\wedge\grw_{j})\|_{\grr_{j+1},0} & \leq &
   \gry_{j}^{(1)}\|\grw_{j}\wedge\grw_{j}\|_{\grr_{j},1}
       \leq  2 \tilde{c}_{0}\gry_{j}^{(1)}\|\grw_{j}\|_{\grr_{j},0}
        \|\grw_{j}\|_{\grr_{j},1}.  \nonumber
\end{eqnarray}
We add together all first order derivatives $\del_{x}$. Then (3.14)
with $k = 0$, gives
\begin{equation}
   \grd_{j+1}^{(1)} \leq 2\tilde{c}_{0}\{2\tilde{c}_{0}\gry_{j}^{(1)}
   \grd_{j}^{(0)} + \tilde{c}_{0}\gry_{j}^{(1)}(\grd_{j}^{(0)} +
   \gry_{j}^{(0)}(\grd_{j}^{(0)})^{2}) + \tilde{c}_{0}\gry_{j}^{(0)}
   \grd_{j}^{(0)}\}\grd_{j}^{(1)} + \gry_{j}^{(0)}(\grd_{j}^{(0)})^{2}.
\end{equation}
We have
\begin{equation}
  \gry_{j}^{(k)}/\gry_{j}^{(k-1)} \leq \hat{c}_{k}4^{j},
\end{equation}
where $\hat{c}_{k} = 4c_{k+3}/(c_{k+2}\grr_{\infty}^{2})$. Taking
$k=1$ in this gives
\begin{equation}
    \grd_{j+1}^{(1)} \leq \grg_{j}^{(1)}\grd_{j}^{(1)}, \; \;
    \grg_{j}^{(1)} = 2\tilde c_0(4\tilde{c}_{0}\hat{c}_{1}4^{j} + \tilde{c}_{0} + 1)
    \gry_{j}^{(0)}\grd_{j}^{(0)}.
\end{equation}
Eventually $\grg_{j}^{(1)}\rightarrow 0$ rapidly.  It follows that
from some $j_{1}$ onward, $\grd_{j}^{(1)}$, and
$\gry_{j}^{(1)}\grd_{j}^{(1)}$ decrease rapidly to zero, and we have
the results  of the previous section for $k=1$.

We repeat this process.  We take another first order derivative in
(4.1). Similar arguments lead to
\begin{equation}
     \grd_{j+1}^{(2)} \leq \grg_{j}^{(2)}\grd_{j}^{(2)},
\end{equation}
where eventually $\grg_{j}^{(2)}\rightarrow 0$ rapidly, and
convergence is established with $k=2$, and so on.  Thus, we see that
$\grw_{j}$ and $B_{j}$ tend to zero rapidly with all derivatives up
to order $k$, and that $G_{\infty}\in C^{k}(M_{\infty})$. This
proves theorem (0.1), except that we only have $A=G_{\infty}\in
C^{k}$.

\section{H\"older continuity}
\setcounter{equation}{0}

We now assume that $M$ is of class $C^{l}$, $l\geq4$, and that our
original $\grw=\grw_{0}$ is of class $C^{k}$, $1\leq k\leq l-4$, and
the constructed sequences $B_{j}$, $\grw_{j}$ are converging rapidly
to zero in $C^{k}$-norm.  We want to show that the $B_{j}$ converge
rapidly to zero in the standard H\"older $C^{k,\gra}$-norm,
$0\leq\gra\leq1/2$,
\begin{equation}
  \|B\|_{C^{k,\gra}(D)} = \|B\|_{C^{k}(D)} + H_{\gra,D}(\del^{k}B).
\end{equation}
Here $\del^{k}B$ stands for all $k$-th order derivatives of $B$.

This follows, in principle, from the well known $1/2$-estimate [4],
applied to the operators $P_{0}$, $P_{0}^{(K)}$ in (2.2), acting on
forms with compact support. (A precise version will be given in
[3].) Thus we shall assume the following, where the set-up is as in
section 2.
\begin{lemma}
  For all $\grph\in C^{k}(D)$ of compact support, we have
\[
   \|P_{0}\grph\|_{C^{k,\gra}(D)} \leq c_{k+4}\|\grph\|_{C^{k}(D)},
\]
where $0\leq\gra\leq1/2$, and $k\geq0$ is an integer.
\end{lemma}

As in section 2 we use the cutoff function $\grl_{0}$. We sacrifice
$1/2$ derivative for the sake of simplicity to get
\begin{equation}
  \|P_{0}\grph\|_{C^{k,\gra}(D_{0})} \leq c_{k+4}\|\grl_{0}\grph\|_{C^{k}(D)}
   + \|P_{0}(\grl_{1}\grph)\|_{C^{k+1}(D_{0})}.
\end{equation}
Also
\begin{equation}
  \|P_{1}\grph\|_{C^{k,\gra}(D_{0})} \leq \|P_{1}\grph\|_{C^{k+1}(D_{0})}.
\end{equation}
This leads to
\begin{equation}
   \|P\grph\|_{\grr(1-\grs),k,\gra} \leq \gry\|\grph\|_{\grr,k},
\end{equation}
for $0\leq\gra\leq 1/2$. Replacing $k$ by $k+1$ in (2.8) gives
\begin{equation}
   \gry = c_{k+4}\grs^{-2n-2(k+1)+1}\grr^{-2(k+1)}.
\end{equation}

With these estimates we readily see that the sequence $B_{j}$
rapidly decreases in $C^{k,\gra}(M_{\infty})$ in the previous
argument, and that the limit $A=G_{\infty}\in
C^{k,\gra}(M_{\infty})$, $0\leq\gra\leq1/2$. This proves theorem
(0.1), in full.


\section{Scale invariance on the Heisenberg group}
\setcounter{equation}{0}

The difficulties of the last two sections stem mainly from the blow
up of the coefficients in our estimates as we initially shrink the
domain.  This can be largely overcome on the Heisenberg group, or
real hyperquadric, (1.1) with $h=0$,
\begin{equation}
M : y^{n}=|z'|^{2}, \; \; \;  M_{\grr} = M \cap \{|z^{n}|\leq\grr
\},
\end{equation}
by using scale invariance in two ways.  Each $M_{\grr}$ is a graph
over the corresponding Heisenberg ball $D_{\grr}$, $|z'|^{4} +
(x^{n})^{2}\leq\grr^{2}$.  These domains are permuted by the
non-isotropic dilations [2],
\begin{equation}
  T_{\grr}(z',z^{n}) = (\sqrt{\grr}z', \grr z^{n}), \; \; \;
  M_{\grr}=T_{\grr}(M_{1}),
\end{equation}
which also preserve the CR structure.

\vspace{2ex} \noindent\textbf{A)} By the explicit form of their
kernels as given in [11], the operators $P_{0}$, $Q_{0}$, $P_{1}$,
$Q_{1}$ are easily seen to be invariant under these scalings,
\begin{equation}
   P_{(\grr)}= T^{*}_{\grr^{-1}}P_{(1)}T^{*}_{\grr}, \; \;
   Q_{(\grr)}= T^{*}_{\grr^{-1}}Q_{(1)}T^{*}_{\grr}.
\end{equation}
In fact, (modifying somewhat  the notations of (1.30) [11]),
change-of-variables gives
\begin{eqnarray}
  (P_{0(\grr)}\grph)(z) & = & \int_{M_{\grr}}\grph(\grz)\wedge\grW(\grz,z)
     =  \int_{M_{1}}(T_{\grr}^{\grz})^{*}(\grph\wedge\grW)   \\
    & = &
  (T_{\grr^{-1}}^{z})^{*}\int_{M_{1}}(T_{\grr}^{\grz})^{*}(\grph)\wedge
  (T_{\grr}^{\grz}\times T_{\grr}^{z})^{*}\grW.  \nonumber
\end{eqnarray}
But $(T_{\grr}^{\grz}\times T_{\grr}^{z})^{*}\grW = \grW$, which
also holds for $\grW_{1}$, the kernel (1.31) [11] of the boundary
integral $P_{1}$.

We define scale-invariant Folland-Stein H\"older norms by
\begin{equation}
   \|\grph\|_{\grr,k,\gra} = \|T^{*}_{\grr}\grph\|_{1,k,\gra},
\end{equation}
where the right-hand side is the ordinary FS norm on $M_{1}$. As
soon as we have fixed a positive lower bound on $\grr$,  these norms
are equivalent to the ordinary Folland-Stein norms [2].

As noted in [4], on the Heisenberg group the operators $P_{0}$,
$Q_{0}$ are essentially divergences of the Folland-Stein fundamental
solution operators.  In fact, their kernels can be constructed from
the functions (6.1) in [2] by taking one Folland -Stein derivative.
We may appeal to theorem (10.1) of [2].  For compactly supported
forms $\grph$ on the unit Heisenberg ball $M_{1}$, and $k\in\bZ$,
$0<\gra<1$, this gives
\begin{equation}
   \|P_{0(1)}\|_{1,k+1,\gra} \leq c_{k,\gra}\|\grph\|_{1,k,\gra},
\end{equation}
and similarly for $Q_{0(1)}$.

We then follow the procedure of section 2 using a cutoff function
and crudely estimating the boundary integrals, but using
Folland-Stein norms, on the unit Heisenberg ball $M_{1}$. We take
one extra derivative in (2.5) and (2.6) for the H\"older ratio.  The
denominators in the coefficients for the estimates over $M_{1-\grs}$
involve only powers of $\grs$, and this remains so after scaling.
Thus we get
\begin{equation}
   \|P_{(\grr)}\|_{\grr(1-\grs),k+1,\gra} \leq \gry\|\grph\|_{\grr,k,\gra},
\end{equation}
\begin{equation}
  \gry = c_{k,\gra}\grs^{-s(k,n)},
\end{equation}
where $s(k,n)$ is a positive integer.

To prove theorem (0.2) we apply the procedure of section 3 with the
above constructions.  We see that we gain one Folland-Stein
derivative in passing from $\grw$ to $B$.  Initially we normalize
$\grw_{0}$ as in lemma (1.2) to order $k$.  Then
$\|\grw_{0}\|_{\grr_{0},k,\gra}$, and hence
$\|B_{0}\|_{\grr_{0},k+1,\gra}$ tend to zero as we shrink
$\grr_{0}$.  Thus we need no extra smoothness to get the process of
section 3 started. It yields theorem (0.2) in one step.

Strictly speaking, we need  $k\geq 2$ so that our forms are at least
of class $C^{1}$ for the direct derivation of the homotopy formula
in [11]. However, on the Heisenberg group we can establish the
homotopy formula directly from the Folland-Stein constructions [2].
This only requires $k\geq 1$.

\vspace{2ex} \noindent\textbf{B)} We indicate another method using
the ordinary norms of sections 2 and 5. This uses the dilations
$T_{\grk}$, $\grk > 0$, to pull back the (trivial) vector bundle and
the connection from $M_{\grk}$ to $M_{1}$,
\begin{equation}
  \grw^{\grk}\equiv T_{\grk}^{*}\grw = \sqrt{\grk}\sum_{\gra=1}^{n-1}
   \grG_{\barr{\gra}}(\sqrt{\grk}z',\grk x^{n}) d\barr{z}^{\gra}.
\end{equation}
The coefficients, together with any $(z',x^{n})$-derivatives tend to
zero uniformly on $D_{1}$, as $\grk\rightarrow 0$.

We set $\grr_{0}=1$ in (3.1), getting a fixed decreasing sequence of
radii $\grr_{j}\rightarrow\grr_{\infty}>0$, and the corresponding
sequence of Heisenberg balls $D_{\grr_{j}}$, $0\leq j\leq\infty$. We
also choose cut-off functions $\grl_{j}=\grl_{j}(|z^{n}|)\geq0$,
$\grl_{j}=1$ on $D_{\grr_{j}}$, $\grl_{j}=0$ near $\del
D_{\grr_{j-1}}$. We have $|\del^{(1)}\grl_{j}|\leq
c/(\grr_{j}\grs_{j})\leq \tilde{c}/\grs_{j}$, where
$\tilde{c}=c/\grr_{\infty}$, etc.

Since the geometry of Heisenberg balls is different from that of
Euclidean balls, the corresponding quantities $\grd(D)$,
$\grd(D_{0},\del D)$, and vol$(D-D_{0})$ of section 2 are different.
For each $j$, they depend on $\grr_{j}$, $\grs_{j}$.  But since we
have the fixed lower bound $\grr_{j}\geq\grr_{\infty}$, we can again
replace (2.8), (5.5) by
\begin{equation}
  \gry_{k} = c_{k+4}\grs^{-s(n,k)},
\end{equation}
where $s(n,k)$ is a positive integer.

Now we take $\grw_{0}=\grw^{\grk}$, so that by (6.9)
$\grz_{0}=O(\sqrt{\grk})$ in (3.16).  By taking $\grk>0$
sufficiently small, we can start the argument of sections 3 and 5.
It yields convergence on $D_{\infty}$.  Transforming back by
$T_{\grk}^{-1}$ gives theorem (0.1) directly (for $k$ finite),
without recourse to sections one or four.


\bigskip
\noindent Xianghong Gong \\
    Department of Mathematics\\University of Wisconsin, Madison, WI
    53706\\
   E-mail: gong@math.wisc.edu \\

   \bigskip
\noindent S. M. Webster \\
    Department of Mathematics\\University of Chicago,  Chicago, IL
    60637\\
   E-mail: webster@math.uchicago.edu \\


\begin{thebibliography}{99}

\bibitem{1} S-C. Chen and M-C. Shaw, Partial Differential
Equations in Several Complex Variables, AMS/IP, Studies in Adv.
Math. 19 (2001).

\bibitem{2} G. B. Folland and E. M. Stein, Estimates for the $\barr{\del}_{b}$
complex and analysis on the Heisenberg group, Comm. Pure Appl. Math.
27 (1974) 429-522.

\bibitem{3} X. Gong and S. M. Webster, Regularity for the CR vector bundle
problem II, pre-print.

\bibitem{4} G. M. Henkin, The Lewy equation and analysis on pseudoconvex
domains, Russ. Math. Surv. 32 (1977) 59-130.

\bibitem{5} J-L. Koszul and B. Malgrange, Sur certaines structures fibr\'ees
complexes, Arch. Math. 9 (1958) 102-109.

\bibitem{6} L. Ma and J. Michel, Regularity of local embeddings of strictly
pseudoconvex CR structures, J. Reine Angew. Math. 447 (1994)
147-164.

\bibitem{7} L. Ma and J. Michel, On the regularity of CR structures for
almost CR vector bundles, Math. Z. 218 (1995) 135-142.

\bibitem{8} J. K. Moser, A rapidly convergent iteration method and nonlinear
differential equations I, Ann. Scuola Norm. Pisa 20 (1966) 265-315.

\bibitem{9} A. Nagel and J-P. Rosay, Non existence of a homotopy formula
for (0,1) forms on hypersurfaces in $\bC^{3}$, Duke Math. Jour. 58
(1989) 823-827.

\bibitem{10} C. Romero, Potential theory for the Kohn Laplacian on the
Heisenberg group, thesis, University of Minnesota (1991).

\bibitem{11} S. M. Webster, On the local solution of the tangential
Cauchy-Riemann equations, Ann. Inst. H. Poincare 6 (1989) 167-182.

\bibitem{12} S. M. Webster, On the proof of Kuranishi's embedding theorem,
Ann. Inst. H. Poincare 6 (1989) 183-207.

\bibitem{13} S. M. Webster, The integrability problem for CR vector bundles,
Proc. Symp. Pure Math. 52 (1991) part 3, 355-368.




\end{thebibliography}
\end{document}